\documentclass{article}

\usepackage{amssymb,amscd,amsthm,latexsym}
\usepackage[all]{xy}

\newtheorem{defi}{Definition}[section]
\newtheorem{theo}{Theorem}[section]
\newtheorem{lemma}{Lemma}[section]
\newtheorem{prop}{Proposition}[section]
\newtheorem{coro}{Corollary}[section]

\def\into{ \rightarrowtail }
\def\onto{ \twoheadrightarrow }
\def\splito{ \rightleftarrows }
\newdir{ >}{{}*!/-8pt/@{>}}
\def\EE{ \mathbb{E} }
\def\CC{ \mathbb{C} }
\def\DD{ \mathbb{D} }
\def\VV{ \mathbb{V} }

\def\NN{ \mathbb{N} }
\def\TT{ \mathbb{T} }

\begin{document}

\author{Dominique Bourn}

\title{On the cocartesian image of preorders and equivalence relations in regular categories}

\date{}

\maketitle

\begin{abstract}
In a regular category $\EE$, the direct image along a regular epimorphism $f$ of a preorder is  not a preorder in general. In 
$Set$, its best preorder approximation is then its cocartesian image above $f$. In a regular category, the existence of such a cocartesian  image above  $f$ of a preorder $S$ is actually equivalent to the existence of the supremum $R[f]\vee S$ among the preorders. We investigate here some conditions ensuring the existence of these cocartesian images  or equivalently of these suprema. They applied to two very dissimilar contexts: any topos $\EE$ with suprema of chains of subobjects or any $n$-permutable regular category.	
\end{abstract}

\section*{Introduction}

In a regular category $\EE$, the direct image of a reflexive (resp. symmetric) relation along a regular epimorphism $f$ is a reflexive (resp. symmetric) one. However, this is no more the case, in general, for preorders and a fortiori for equivalence relations. The best preorder approximations are then their cocartesian images above $f$ in the respective categories $PrO\EE$ of preorders and $Equ\EE$ of equivalence relations in $\EE$.

The following observation was made in \cite{BGour}, emphasizing the importance of this notion, but was not much developped, the core of that work being elsewhere:\\
\noindent\textbf{Proposition.} \emph{Let $\EE$ be a regular category and $f:X\onto Y$ a regular epimorphism.  The following conditions are equivalent:\\
	1) the map $\phi: R\to T$ is cocartesian in $Equ\EE$ above  $f$ in $\EE$;\\
	2) the map $\phi: R\to T$ is a regular epimorphism in $Equ\EE$ above  $f$ in $\EE$;\\
    3) we have: $f^{-1}(T)=R[f]\vee R$ in $Equ\EE$.}

In particular, the question of the conditions ensuring the existence of such cocartesian maps (or equivalently the  existence of suprema of equivalence relations) was not addressed. And of course, it is an important question, since it would structurally allow to extend the equivalence relations along any regular epimorphism. The purpose of this work is to investigate it. And since the hereabove  result is easily shown to remain valid when $R$ is only a  preorder, we  shall extend our  investigation to cocartesian images (or suprema of pairs) of preorders. 

In \cite{CKP}, the existence of the supremum  $R\vee S$ of two equivalence relations in a regular category $\EE$ is strongly related to the following chain of reflexive relations:
$$((R,S)):  \Delta_X \subset R \subset RS=(R,S)_2 \subset RSR=(R,S) _3\subset RSRS=(R,S)_4 \subset \cdots  $$
i) First, it is noticed that, in the case of varieties, the supremum $R\vee S$ among the equivalence relation is the union of this chain. So, according to the previous proposition, any variety has cocartesian images.\\
ii) It is noticed as well (Theorem 3.1), that, in a regular category $\EE$, for two equivalence relations $R$ and $S$, the following two conditions are equivalent:\\
\emph{1) the chain $((R,S))$ is stationary from level $n$;\\
2) the reflexive relation $(R,S)_n$ is an equivalence relation.\\
Under any of these conditions, we get: $(R,S)_n=R\vee S$.}\\
We shall start from these two observations, however 
extended from equivalence relations to reflexive relations, to produce the conditions we are looking for. 

Inspired by i), we shall show that, in any regular category $\EE$ in which the chain $((R,S))$, where the pair $(R,S)$ is only a pair of reflexive  relations, has a supremum $\Sigma(R,S)$ among the  reflexive relations, this reflexive relation is actually a preorder, provided that this supremum  is stable under composition of reflexive relations. Accordingly, under this condition, 1) we get $\Sigma(R,S)=R\vee S$ among the preorders, and 2) the inclusions $PrO\EE \hookrightarrow Ref\EE$ and $Equ\EE \hookrightarrow Ref\EE$ have left adjoints, respectively given by $T\mapsto \Sigma(T,T)$ and $T\mapsto \Sigma(T,T^{op})$. This point 2) produces, in $\EE$, a clear construction of the peorder and the equivalence relation generated by any  reflexive relation $T$.

Any elementary topos $\EE$ satisfies this condition provided that it has suprema of chains of subobjects; this is obviously the case of any Grothendieck topos. Moreover, this condition can be then transfered to any category $\TT(\EE)$ of  internal $\TT$-algebras in $\EE$, where $\TT$ is the theory of a finitary algebraic structure. We incidentally show that, if the associated variety $\VV(\TT)$ is congruence modular, so is the category $\TT(\EE)$. 

Inspired by ii), we shall investigate the case in which the chain $((R,S))$, again when the pair $(R,S)$ is only a pair of reflexive  relations, is stationary from level $n$, which obviously implies that $\Sigma((R,S))=(R,S)_n$. We shall show that, in any regular category $\EE$, given any pair $(R,S)$ of reflexive relations on $X$ and $2\leq n$, the following conditions are equivalent:\\ 
\emph{1) $(R,S)_{n+2}=(R,S)_{n}$;\\
2) the chain $((R,S))$ is stationary from $n$;\\
3) the reflexive relation $(R,S)_{n}$ is a preorder.}\\
So, under any of these conditions, when $R$ and $S$ are preorders, this $(R,S)_n$ becomes the supremum $R\vee S$ among the preorders. Accordingly, when $S$ is a preorder such that $(R[f],S)_{n+2}=(R[f],S)_{n}$, the existence of the cocartesian image of $S$ along the regular epimorphism $f:X\onto Y$ is ensured. This implies, in particular, that any congruence $n$-permutable regular category (namely any regular categories in which $(R,S)_n=(S,R)_n$ for any pair $(R,S)$ of equivalence relations) has cartesian images of equivalence relations along regular epimorphisms.
As  a collateral effect, will unexpectedly follow a characterization of congruence $n$-permutable regular categories in terms of existence of cocartesian images which is valid when $n$ is odd and which only gives rise to an implication when $n$ is even: there lies a surprising asymmetry between the odd and even cases.

The article is organized along  the following lines. Section 1 is devoted to basic recalls on reflexive relations. Section 2 investigates the distinction between direct image and cocartesian image. Section 3 investigates the properties of the chain $((R,S))$ where $R$ and $S$ are only reflexive relations and in particular the conditions under which it admits a supremum among the reflexive relations. Section 4 is devoted to showing that any elementary topos $\EE$ admits cocartesian images of preorders  and equivalence relations provided that any chain of subobjects has a supremum, and Section 5  to producing a direct application concerning a transfer of congruence modularity. Section 6 is devoted to showing that any $n$-permutable regular category admits cocartesian images of preorders and equivalence relations, and this gives rise to three different kinds of characterization of these categories: Theorem \ref{disy}, Corollary \ref{genmal} and Theorem \ref{mixed}.

\section{Quick recalls on reflexive relations}

In this article we shall assume that any category $\EE$ is finitely complete. 

\subsection{Basics about relations}

As usual, a relation $R$ between two objects $X$ and $Y$ in a category $\EE$ is a subobject $(d_0^R,d_1^R): R\into X\times Y$, or more precisely a class up to isomorphism of this kind of subobjects. Clearly the set $Rel\EE(X,Y)$ of such relations is preordered by the inclusion. The dual of a relation $R$, denoted by $R^{op}$, is given by the subobject $(d_1^R,d_0^R): R\into Y\times X$. As soon as, in addition, the category $\EE$ is regular \cite{Barr}, the relations can be composed:
given any other relation $(d_0^S,d_1^S): S\into Y\times Z$, the composition $S\circ R$ is classically (see \cite{CKP} for instance) given by the following construction:\\
1) first take the pullback $R\times_Y S$ of $d_0^S$ along $d_1^R$:
$$ \xymatrix@=20pt{
	R\times_Y S \ar@<2ex>[rrd]^{p_1^S} \ar@{->>}[rd]^{\psi} \ar@<-2ex>[ddr]_{p_0^R}\\
	& S\circ R \ar@{.>}@(u,)[rr] \ar@{.>}@(l,)[dd]^{} & S \ar[r]_{d_1^S}  \ar[d]^{d_0^S} & Z\\
    &	R \ar[r]_{d_1^R}  \ar[d]^{d_0^R} & Y   &\\
	& X  &
}
$$
2) then take the canonical decomposition into regular epimorphism and mono- morphism of the map $(d_0^R.p_0^R,d_1^S.p_1^S): R\times_Y S \onto S\circ R \into X\times Z$. We then get $(S\circ R)^{op}=R^{op}\circ S^{op}$. This composition is associative up to isomorphism, and makes the discrete (equivalence) relations $\Delta_X$ ($=(1_X,1_X): X\into X\times X$) on $X$ the units of this law. So that the composition $\circ$ produces a category denoted by $Rel\EE$. Compositions on left and right preserve the inclusion.

\medskip

Given any object $X$, we denote by $\nabla_X$ the undiscrete equivalence relation $1_{X\times X}: X\times X \into X\times X$; it is the largest relation on $X$. 
As we know, a relation $R$ on $X$ is \emph{reflexive} if and only if $\Delta_X \subset R$, \emph{symmetric} if and only if $R^{op}\subset R$, and \emph{transitive} if and only if $R^2=R\circ R\subset R$. An internal \emph{preorder} is a reflexive and transitive relation; it is an internal equivalence relation, when, moreover, it is symmetric. For any reflexive relation $S$ on $X$, we chose the simplicial notations:
$$\xymatrix{ S \ar@<-6pt>[r]_{d_1^S} \ar@<6pt>[r]^{d_0^S} & X \ar[l]|{s_0^S} }$$
Given any map $f: X\to Y$ in $\EE$, we define the inverse image $f^{-1}(S)$ of a relation $S$ on $Y$ by the following left hand side pullback in $\EE$:
$$ \xymatrix@=20pt{
	f^{-1}(S) \ar[r]^{\check f} \ar@{ >->}[d]_{(d_0,d_1)} & S\ar@{ >->}[d]^{(d_0^S,d_1^S)} && f^{-1}(S) \ar[r]^{\check f} \ar@{ >->}[d]_{(d_0,d_1)} & S\ar@{ >->}[d]^{(d_0^S,d_1^S)}\\
	X\times X \ar[r]_{f\times f}  & Y\times Y && 	\nabla_X \ar[r]_{\nabla_f}  & 	\nabla_Y
}
$$
The inverse image $R[f]=f^{-1}(\Delta_Y)$ is called the \emph{kernel equivalence relation} of the morphism $f$. The inverse images preserve the reflexive, the symmetric and the transitive relations. They preserve the intersection: $f^{-1}(S)\wedge f^{-1}(R)=f^{-1}(S\wedge R)$ as well, and thus the inclusion. When, moreover,  $\EE$ is regular, we get: $f^{-1}(S)\circ f^{-1}(R) \subset f^{-1}(S\circ R)$.

Given any map $f:X\to Y$, and a pair $(R,S)$ of relations on $X$ and $Y$ respectively, we say that $f$ produces a morphism of relations $R\to S$ above $f$ when $R\subset f^{-1}(S)$. We denote by $Ref\EE$ the induced category which is finitely complete and by $(\;)_0: Ref\EE \to \EE$ the left exact functor associating with any reflexive relation $R$ on $X$ its ground object $X$. It is a fibration whose \emph{cartesian maps} are given by the inverse images. Now, the inverse image $f^{-1}(S)$ can be obtained by the hereabove right hand side pullback in the category $Ref\EE$. Since any of its fibers is a preorder, the functor $(\;)_0:$ is faithful, in other words, between two reflexive relations, there is atmost one morphism above a given morphism $f$ in the category $\EE$.

\medskip

We denote by $PrO\EE$, $RSym\EE$, and $Equ\EE$ the respective subcategories of internal preorders, reflexive and symmetric relations, and equivalence relations. The restrictions of the functor $(\;)_0$ to these three subcategories determine sub-fibrations. Very important and straightforward is the following:
\begin{theo}\label{main3}
	In any category $\EE$, given a split epimorphism $(f,s):X\splito Y$, the inverse image $f^{-1}:PrO_Y\mathbb E\to PrO_X\mathbb E$ (resp. $f^{-1}:Equ_Y\mathbb E\to Equ_X\mathbb E$) induces a preorder preserving bijection between the preorders (resp. the equivalence relations) on $Y$ and the preorders (resp. the equivalence relations) on $X$ containing $R[f]$. The inverse mappings are given by the restriction of $s^{-1}$.
	
	When, moreover, the category $\EE$ is regular we get $f^{-1}f(T)=R[f]\circ T\circ R[f]$; so, the previous bijections can be extended to any regular epimorphism $f:X\onto Y$, the inverse mapping being given by the restriction of the direct image.
\end{theo}
\proof
Thanks to the Yoneda embedding, it is enough to check the first part of the theorem in the category $Set$ of sets. When $T$ is transitive and $R[f]\subset T$, we get:
$$xTx'\iff sf(x)Tsf(x') \iff xf^{-1}(s^{-1}(T))x'$$
As for the second part of the theorem, see the next section for the definition of the direct image. Again a proof in $Set$ is enough thanks to the embedding theorem for regular categories \cite{Barr}. For any reflexive  relation $T$ on $X$, we get:
$$xf^{-1}(f(T))x'\iff f(x)f(T)f(x') \iff x(R[f]\circ T\circ R[f])x'$$
The equivalence $x(R[f]\circ T\circ R[f])x'\iff xTx'$ holds whenever $R[f]\subset T$ and $T$ is transitive.
\endproof

\subsection{The square construction}

Given any pair $(R,S)$ of reflexive relations on an object $X$ in a  category $\mathbb E$, we denote by $R\square S$ the inverse image of the equivalence relation $S\times S$ on $X\times X$ along the inclusion $(d_0^R,d_1^R):R\rightarrowtail X\times X$, see \cite{BGGGum} for instance. In set-theoretical terms, this double relation $R \square S$ is the subset of elements $(u,v,u',v')$ of $X^4$ such that the following relations $uRu',vRv',uSv,u'Sv'$ hold:
$$ \xymatrix@=15pt{
	u \ar@{.>}[r]^S \ar@{.>}[d]_R & v\ar@{.>}[d]^R\\
	u' \ar@{.>}[r]_S & v'
}
$$
This \emph{square construction} defines a double reflexive relation in $\EE$:
$$ \xymatrix@C=3pc@R=2pc{
	R \square S \ar@<-1,ex>[d]_{\delta_0^R}\ar@<+1,ex>[d]^{\delta_1^R} \ar@<-1,ex>[r]_{\delta_0^S}\ar@<+1,ex>[r]^{\delta_1^S}
	& S \ar@<-1,ex>[d]_{d_0^S}\ar@<+1,ex>[d]^{d_1^S} \ar[l]\\
	R \ar@<-1,ex>[r]_{d_0^R} \ar@<+1,ex>[r]^{d_1^R} \ar[u]_{} & X\ar[u]_{} \ar[l]
}
$$
which is actually the largest double reflexive relation on $R$ and $S$. In $Ref\EE$, we get: $R \square S=(d_0^R)^{-1}(S) \cap (d_1R)^{-1}(S)$, and the diagram can be seen as an internal reflexive relation in $Ref\EE$ above the reflexive relation $R$ in $\EE$. The following observations are straightforward:
\begin{lemma}\label{square1}
	Given any map $f:X\to Y$ in $\mathbb E$, and any reflexive relation $T$ on $X$, the kernel equivalence relation in the category $Ref\EE$ of the following map: $T\to \nabla_Y$ above $f$ is given by the following diagram:
	$$ \xymatrix@=35pt{
		R[f] \square T  \ar@<-1,ex>[r]_{\delta_0^T}\ar@<+1,ex>[r]^{\delta_1^T} & T \ar[rr]^{(f,(f.d_0^T,f.d_1^T))} \ar[l] && \nabla_Y
	}
	$$
	Moreover it is the kernel equivalence relation of any other map with domain $T$ above $f$ in $Ref\EE$. Accordingly this equivalence relation is the unique effective equivalence relation on $T$ in $Ref\EE$ above the equivalence relation $R[f]$.
\end{lemma}
\begin{lemma}\label{square2}
Given any pair $(R,T)$ of an equivalence relation and a preorder in $\EE$, the following conditions are equivalent:\\
1) $R\subset T$;\\
2) the following diagram in $Ref\EE$ is made of inverse images (=cartesian maps):
$$ \xymatrix@C=3pc@R=2pc{
	R \square T \ar@<-1,ex>[r]_{(d_0^R,\delta_0^T)}\ar@<+1,ex>[r]^{(d_1^R,\delta_1^T)}
	& T \ar[l]
}
$$
\end{lemma}
\proof
The implication  $1) \Rightarrow 2)$ comes from the fact that when $R$ is an equivalence relation and $T$ a preorder, if we have $R\subset T$, then, in the following situation:\\
$$ \xymatrix@=15pt{
	u \ar@{.>}[d]_R & v\ar@{.>}[d]^R\\
	u'  & v'
}
$$
we get: $uTv\iff u'Tv'$. Conversely, when $(d_0^R,\delta_0^T)$ is an inverve image, the following diagram:
$$ \xymatrix@=15pt{
	u \ar@{.>}[r]^T \ar@{.>}[d]_R & u\ar@{.>}[d]^R\\
	u  & v
}
$$
shows that $uRv\Rightarrow uTv$.
\endproof

\section{Direct image and cocartesian image}

\subsection{Direct image along regular epimorphisms}

When $\EE$ is a regular category, $f: X\onto Y$ a regular epimorphism in $\EE$, and $R$ a relation on $X$, the direct image $f(R)$ of $R$ along $f$ is the relation on $Y$ given by the following canonical decomposition of $(f\times f).(d_0^R,d_1^R): R\to Y\times Y$:
$$ \xymatrix@=25pt{
{R\;}\ar@{ >->}[d]_{(d_0^R,d_1^R)} \ar@{->>}[rr]^{\hat f}	&& f(R)  \ar@{ >->}[d]^{(\delta_0,\delta_1)} \\
	 X\times X \ar@{->>}[rr]_{f\times f} && Y\times Y
}
$$
The direct images preserve the reflexive and the symmetric relations, but not the transitive one. The regular categories in which the direct images preserves the equivalence relations are characterized in \cite{GR}:
\begin{prop}\label{n=1}
When the category $\EE$ is a regular category, the following conditions are equivalent:\\
1) 	the direct images preserves the equivalence relations;\\
2) $\EE$ is a Goursat category, i.e. a category in which any pair $(R,S)$ of equivalence relations on an object $X$ is congruence $3$-permutable, namely satisfies: $S\circ R\circ S=R\circ S \circ R$.
\end{prop}
Straightforward is the following:
\begin{lemma}\label{lemma1}
	In a category $\EE$, given any equivalence relation on $X$, the direct image of $R[d_0^R]$ along the split epimorphism $d_1^R: R\to X$ is nothing but $R$ itself. 
	
	When $\EE$ is regular and $S$ any reflexive relation on $X$, we get $d_1^S(R[d_0^S])=S\circ S^{op}$. Accordingly, in this context, given any reflexive  relation  $R$, it is an equivalence relation if and only if $d_1^R(R[d_0^R])=R$. 
\end{lemma}
\begin{prop}\label{idcomp}
	Let $\EE$ be a regular category, $f:X\onto Y$ a regular epimorphism. 
	Given any pair $(R,S)$ of reflexive relations on $X$, we get:
	$$f^{-1}(f(S)\circ f(R))=R[f]\circ S\circ R[f]\circ R \circ R[f] \;\; (*)$$
	 Whence: $f(S\circ R)\subset f(S)\circ f(R)= f(S\circ R[f]\circ  R)$. Accordingly we get: $f(R)^{n+1}=f(R\circ R[f]\circ R \cdots \circ R[f]\circ R)$, with $n$ times $R[f]$.
\end{prop}
\proof Thanks to the embedding theorem for regular  categories \cite{Barr}, it is enough to check all this of in the category $Set$ of sets, where the proofs of the first assertion is straightforward. The direct image along $f$ of the inclusion:\\ $S\circ R\subset R[f]\circ S\circ R[f]\circ R \circ R[f]$ produces: $f(S\circ R)\subset f(S)\circ f(R))$. 

Now, since $f$ is a regular epimorphism, we get $f(S)\circ f(R)= f(S\circ R[f]\circ  R)$ if and  only if: $f^{-1}(f(S)\circ f(R))= f^{-1}(f(S\circ R[f]\circ  R))$. We know that the first term is  $R[f]\circ S\circ R[f]\circ R \circ R[f]$ by $(*)$, while the second one is $R[f]\circ (S\circ R[f]\circ R) \circ R[f]$ by Theorem \ref{main3}.
\endproof

\subsection{Cocartesian image along regular epimorphisms}

Now let us recall the following:

\begin{defi}
Given any functor $U: \CC \to \DD$, a map $f:X\to Z$ in $\CC$ is said to be cocartesian with respect to the functor $U$ when it is universal among the maps with domain $X$ and image $U(f)$, i.e. such that given any other map $\phi:X\to W$ with $U(\phi)=U(f)$, there is a  unique factorization $\psi:Z\to W$ satisfying $\psi.f=\phi$ and $U(\psi)=1_{U(Z)}$. 	
\end{defi}
The fact that the direct images along regular epimorphisms $f: X\onto Y$ preserve the reflexive and the symmetric relations makes cocartesian with respect to $(\;)_0$ the morphism $(f,\hat f):R\onto f(R)$  in $Ref\EE$ and $RSym\EE$ respectively. When $R$ is a preorder (resp. an equivalence relation) and $f(R)$ happens to be so, then this same map is cocartesian in $PrO\EE$ (resp. in $Equ\EE$) as well. So, according to Proposition \ref{n=1} above, in a Goursat regular category $\EE$, the direct and the cocartesain images  along regular epimorphisms do coincide in $Equ\EE$.  

\medskip

When $T$ is a preorder and  the reflexive relation $f(T)$ is not a preorder, arises the question of the existence of a map with domain $T$ which would be cocartesian  with respect to $(\;)_0$ above the regular epimorphism $f: X\onto Y$. When it exists, we shall denote it by $(f,f_!):T\to f_!(T)$ and call its codomain $f_!(T)$ the \emph{cocartesian image} of $T$. The structural importance of this kind of maps is stengthened by the following:
\begin{lemma}
	Let $\EE$ be a regular category and $f: X\onto Y$ a regular epimorphism in $\EE$. The following conditions are equivalent:\\
	1) $(f,\bar f): T\to S$ is cocartesian above $f$ in $PrO\EE$ (resp. in $Equ\EE$);\\
	2) $(f,\bar f): T\to S$ is a regular epimorphism in $PrO\EE$ (resp. in $Equ\EE$).
\end{lemma}
\proof
A regular epimorphism being the cokernel of its kernel equivalence relation, this is a  straightforward consequence of Lemma \ref{square1}.
\endproof

We can extend to preorders a classical observation (see \cite{BGour} for instance) concerning equivalence relations:
\begin{prop}\label{prop10} 
In any category $\EE$, a map $(f,\bar f):T\to S$ is cocartesian in $PrO\EE$ (resp. in $Equ\mathbb E$) above the split epimorphism $(f,s): X\splito Y$ in $\mathbb E$ if and only if we have $f^{-1}(S)=R[f]\vee T$ in $PrO\EE$ (resp. in $Equ\mathbb E$).

When, in addition, the category $\EE$ is regular, this result is valid for any regular epimorphism $f$ as well. In this case we get $f_!(T)=f(R[f]\vee T)$.
\end{prop}
\proof
The fact that a map $(f,\bar f): T\to S$ is cocartesian is characterized by the fact that, for any other morphism $(f,\tilde f):T\to S'$ in $PrO\EE$ (resp. in $Equ\EE$) above $f$, we have $S\subset S'$, or equivalently by the fact that, for any preorder (resp. equivalence relation) $S'$ on $Y$, we have:
$$T\subset f^{-1}(S') \iff S\subset S' \iff f^{-1}(S)\subset f^{-1}(S')$$
the second equivalence being given by Theorem \ref{main3}. Now, the bijection of this same theorem characterizes $T$ by: for any preorder (resp.  equivalence relation) $U$ above $X$ containing $R[f]$, we have: $T\subset U \iff f^{-1}(S)\subset U$ or, in other words, by $f^{-1}(S)=R[f]\vee T$
in $PrO\EE$ (resp. in $Equ\mathbb E$). When $\EE$ is a regular category and $f$ is a regular epimorphism, the same proof remains valid by the second part of this same Theorem \ref{main3}. The last assertion is then straightforward since $R[f]\subset R[f]\vee T$.
\endproof
From that emerges a modular formula for preorders:
\begin{coro}\label{sup}
	Let $\EE$ be a regular category, and $f: X\onto Y$ any regular epimorphism. When $S$ is a preorder on $X$ such that its direct  image $f(S)$ is a preorder on $Y$ as well, we then get: 
	$R[f]\vee S=f^{-1}(f(S))$.\\
	Given any preorder $T$ such  $R[f]\subset  T$, it then follows: $f(S)\wedge f(T)=f(S\wedge T)$. Then the direct image $f(S\wedge T)$ is necessarily a preorder on $Y$,
 and this last equality becomes equivalent to the modular formula for  preorders:
	$(R[f]\vee S) \wedge T=R[f]\vee (S\wedge T)$.
\end{coro}
\proof
Our assumption implies that $(f,\check f): S\onto f(S)$ is cocartesian in $PrO\EE$; whence the first point.

Now, consider any preorder $T$ on $X$ such that $R[f]\subset T$. We know that $f(S\wedge T)\subset f(S)\wedge f(T)$. Now assume that the left hand side diagram is valid:
$$ \xymatrix@=15pt{
	x \ar@{.>}[r]^{R[f]}  \ar@{.>}[d]_S & z  \ar@{.>}[d]^T &&  x \ar@{.>}[d]_S \ar@{.>}@(r,r)[d]^T\\
	x' \ar@{.>}[r]_{R[f]} & z' &&   x'
}
$$
Since $R[f]\subset T$, we get the right hand side situation. This implies $f(S)\wedge f(T) \subset f(S\wedge T)$, whence $f(S\wedge T)=f(S)\wedge f(T)$, and $f(S\wedge T)$ is an equivalence relation since so is $f(S)$.

We then get $f^{-1}(f(S))\wedge T=f^{-1}(f(S))\wedge f^{-1}(f(T))=f^{-1}(f(S)\wedge f(T))=f^{-1}f(S\wedge T)$ which, according to the first part of the proof, is exactly the modular formula.
\endproof
By the previous proposition we obtain a generalization of Propositions 2.1 and 2.4 in \cite{BGour} concerning the supremum of equivalence relations:
\begin{theo}\label{main2}
	For any category $\EE$ the following conditions are equivalent:\\
	1) there are cocartesian maps with respect to $(\;)_0$ in $PrO\EE$  (resp. in $Equ\EE$) above any split epimorphism in $\EE$ and with any domain in $PrO\EE$  (resp. in $Equ\EE$);\\
	2) there are suprema $R\vee S$ in $PrO\EE$  (resp. in $Equ\EE$) of any pair $(R,S)$ of an equivalence relation and a preorder (resp. any pair of equivalence relations). 
	
	When, in addition, the category $\EE$ is  regular, then the following conditions are equivalent:\\
	1) there are cocartesian maps with respect to $(\;)_0$ in $PrO\EE$  (resp. in $Equ\EE$) above any regular epimorphism in $\EE$ and with any domain in $PrO\EE$  (resp. in $Equ\EE$);\\
	2) there are suprema $R\vee S$ in $PrO\EE$  (resp. in $Equ\EE$) of any pair $(R,S)$ of an equivalence relation and a preorder (resp. any pair of equivalence relations).
\end{theo}
\proof
The implication $2) \Rightarrow 1)$ is a consequence of the previous  proposition. Suppose now 1). Starting with a pair $(R,S)$ of an equivalence relation and a preorder, let us show that the codomain  $(d_1^R)_!((d_0^R)^{-1}(S))$ of the cocartesian map in $PrO\EE$ above the split epimorphism $d_1^R:R\splito X$ having  $(d_0^R)^{-1}(S)$ as domain is the supremum  $R\vee T$ in $PrO\EE$. The composition of the maps $(s_0^R, \bar s_0^R): S\into (d_0^R)^{-1}(S)$ and $(d_1^R,(d_1^R)_!):(d_0^R)^{-1}(S) \to (d_1^R)_!((d_0^R)^{-1}(S))$ produces the inclusion $S\subset (d_1^R)_!((d_0^R)^{-1}(S))$. The inverse image along $d_0^R$ of the  inclusion $\Delta_X\subset S$ produces the inclusion $R[d_0^R]\subset (d_0^R)^{-1}(S)$ which in turn produces the inclusion of the cocartesian images $R=(d_1^R)_!(R[d_0^R])\subset (d_1^R)_!((d_0^R)^{-1}(S))$, see Lemma \ref{lemma1}.

Now, suppose that $T$ is a preorder such that $R\subset T$ and $S\subset T$. From $S\subset T$, we first get  $(d_0^R)^{-1}(S)\subset (d_0^R)^{-1}(T)=R\square T$ by Lemma \ref{square2} since $R\subset T$; we then get the inclusion of the cocartesian  image $(d_1^R)_!((d_0^R)^{-1}(S))\subset (d_1^R)_!(R\square T)=T$, this last equality following from the fact that, being split in $PrO\EE$, the morphism $(d_1^R,\delta_1^T):R\square T\to T$ is necessarily cocartesian above its image $d_1^R$ in $\EE$.

When $\EE$ is regular, Condition 1) is valid for any split epimorphism, whence 2), by the first part of the theorem. Suppose 2); take a regular epimorphism $f:X\onto Y$ and $T$ any preorder on $X$, then $f_!(T)=f(R[f]\vee T)$ by Proposition \ref{prop10}.
\endproof

\subsection{A basic example of cocartesian image}\label{specif}

Assume that $R$ is a reflexive relation on $X$, we obtain a chain $((R))$ of reflexive relations on $X$:
$$\Delta_X \subset R \subset R^2 \subset R^3 ... \subset R^n \subset R^{n+1} ...$$ 
When, in addition, $R$ is symmetric, it gives rise to a chain of symmetric relations. A reflexive relation $R$ is transitive if and only if, this chain is stationary from level 1.

\begin{prop}\label{Rn}
	Let $\EE$ be a regular category and $R$ a reflexive (resp. reflexive and symmetric) relation on $X$. The following conditions are equivalent:\\
	1) $R^{n+1}=R^n$,\\
	2) $R^n$ is a preorder (resp. an equivalence relation).
\end{prop}
\proof
1) If we have $R^{n+1}=R^n$, by induction we get $R^m=R^n$ for all $m$ such that $n\leq m$. Now, $R^n\circ R^n=R^{2n}=R^n$ and $R^n$ is a preorder.\\
2) If $R^n$ is a preorder, then $(R^n)^2=R^{2n}\subset R^n$. Then we get: $R^{n+1}\subset R^{2n}\subset R^n$ and $R^{n+1}=R^n$.
\endproof

\begin{coro}\label{firstex}
	Let $\EE$ be a regular category, and $f:X\onto Y$ a regular epimorphism. Let $T$ be a preorder (resp. an equivalence relation) on $X$ such that $f(T)^{n+1}=f(T)^n$. Then the map: $T \stackrel{\hat f}{\onto} f(T)\hookrightarrow f(T)^n$ is cocartesian above $f$ in $PrO\EE$ (resp. in $Equ\EE$). Accordingly we get $R[f]\vee T=f^{-1}(f(T)^n)$ in the respective categories.
\end{coro}
\proof
Let $S$ be a preorder on $Y$ such that $T\subset f^{-1}(S)$. We get:
$f(T)\subset f(f^{-1}(S))=S$ thanks to Theorem \ref{main3} since $S$ is a preorder. And so, we get: $f(T)^n\subset S$.
\endproof

\section{The intertwined chains of pairs of relations}

\subsection{Definition}

From now on, for sake of simplicity, we shall denote by $SR$ the composition  $S\circ R$ of relations. The reflexive relations are stable under composition while the symmetric and the transitive ones are not in general. Now, given any pair $(R,S)$ of reflexive relations on $X$, we shall be interested in the following chain  $((R,S))$ of reflexive relations:
$$\Delta_X \subset R \subset RS \subset RSR \subset RSRS \subset ...$$
which, following the notations of \cite{CKP} for the pairs of equivalence relations, will be, term by term, denoted in the following way:
$$(R,S)_0 \subset (R,S)_1 \subset (R,S)_2 \subset (R,S)_3 \subset (R,S)_4 \subset ...$$
It is clear that the chain $((T))$ of Section \ref{specif} coincides with the chain $((T,T))$. The chain $((R,S))$ and its symmetric one $((S,R))$
are connected by intertwined inclusions:
$$ \xymatrix@=20pt{
((R,S)): & \cdots & (R,S)_n \ar@{ >->}[r]\ar@{ >->}[rd] & (R,S)_{n+1} \ar@{ >->}[r] \ar@{ >->}[rd] & (R,S)_{n+2} & \cdots  \\
((S,R)): & \cdots & (S,R)_n \ar@{ >->}[r] \ar@{ >->}[ru] &  (S,R)_{n+1} \ar@{ >->}[r]\ar@{ >->}[ru] &  (S,R)_{n+2} &  \cdots 
}
$$

Whence for all $1\leq k$: $(S,R)_n \subset (R,S)_{n+k}$ and $(R,S)_{n} \subset (S,R)_{n+k}$ .
Accordingly, if the chain $((R,S))$ is stationary from $n$, $2\leq n$, we get: $(S,R)_n\subset (R,S)_n$ and $(S,R)_{n+k}=(R,S)_{n+k}=(R,S)_n$, for all $1\leq k$. So, the two chains coincide from $n+1$, and, in particular,  the chain $((S,R))$ is stationary from $n+1$. 

It is clear that any preorder $T$ on $X$ such that $R\subset T$ and $S\subset T$ contains $(R,S)_n$ and $(S,R)_n$; so, when $(R,S)$ is a pair of preorders, if one of these two reflexive relations is a preorder, this preorder is certainly the supremum $R\vee S$ among the preorders. Another interest of the previous chains is the following:
\begin{lemma}\label{image}
	Let $\EE$ be a regular category and $f:X\onto Y$ a regular epimorphism. Given any reflexive relation $T$ on $X$, we have:\\
	1) $f((T,R[f])_{2m-1})=f(T)^m=f((R[f],T)_{2m+1})$;\\
	2) $f((R[f],T)_{2m})=f(T)^m=f((T,R[f])_{2m})$.
\end{lemma}
\proof
Given any reflexive relation $T$ on $X$, by Proposition \ref{idcomp}, we have:\\ $f(T)^m=f((T,R[f])_{2m-1})$. On the other hand, we have:\\ $f((R[f],T)_{2m+1})\subset (f(R[f]),f(T))_{2m+1}=(\Delta_Y,f(T))_{2m+1}=f(T)^m$\\ $=f((T,R[f])_{2m-1})\subset f((R[f],T)_{2m+1})$. Whence: $f((R[f],T)_{2m+1})=f(T)^m$.

Similarly we get:\\
1) $f(T)^m=f((T,R[f])_{2m-1})\subset f((R[f],T)_{2m})\subset (\Delta_Y,f(T))_{2m}=f(T)^m$;\\
2) $f(T)^m=f((T,R[f])_{2m-1})\subset f((T,R[f])_{2m})\subset (f(T),\Delta_Y)_{2m}=f(T)^m$.
\endproof
The first observations about the members of this chain are the following ones:
\begin{prop}
	Let $\EE$ be a regular category and $(R,S)$ any pair of reflexive relations on $X$. We get:\\
	1) $(R,S)_{2n}(R,S)_{m}=(R,S)_{2n+m}$; whence $(R,S)_{2n}^k=(R,S)_{2kn}$;\\
	2) $(R,S)_{2(n+m)}\subset (R,S)_{2n+1}(R,S)_{2m}\subset (R,S)_{2n+2m+2}$\\ and $(R,S)_{2n+1}(R,S)_{2m}=(R,S)_{2(n+m)}$ when $R$ is a preorder;\\
	3) $(R,S)_{2(n+m)+1}\subset (R,S)_{2n+1}(R,S)_{2m+1}\subset (R,S)_{2n+2m+3}$,\\ and $(R,S)_{2n+1}(R,S)_{2m+1}=(R,S)_{2(n+m)+1}$, when $R$ is preorder;\\
	4) so, for any pair $(n,m)$ of integers, we get:\\
	$(R,S)_{n+m-1}\subset (R,S)_{n}(R,S)_{m}\subset (R,S)_{n+m+1}$\\
	and for any interger $n$, we get: $(R,S)_{2n-1}\subset (R,S)_n^2\subset (R,S)_{2n+1}$;\\
	5) $(R,S)_{2nk+1}\subset (R,S)_{2n+1}^k\subset (R,S)_{2(n+1)k-1}$ and $(R,S)_{2n+1}^k=(R,S)_{2nk+1}$ when $R$ is a preorder.
\end{prop}
\proof
The first assertion is straightforward. For 2), we check:\\
$(R,S)_{2(n+m)}=(R,S)_{2n}(R,S)_{2m}\subset (R,S)_{2n}R(R,S)_{2m}= (R,S)_{2n+1}(R,S)_{2m}$. When, in addition, $R$ is a preorder, we have $R(R,S)_{2m}=(R,S)_{2m}$; whence the conclusion. On the other hand, we get:\\
$(R,S)_{2n+1}(R,S)_{2m}\subset (R,S)_{2n+1}S(R,S)_{2m}=(R,S)_{2n+2m+2}$.\\
For 3), we check: $(R,S)_{2(n+m)+1}=(R,S)_{2n+1}(S,R)_{2m}\subset (R,S)_{2n+1}R(S,R)_{2m}$ $=(R,S)_{2n+1}(R,S)_{2m+1}\subset(R,S)_{2n+1}S(R,S)_{2m+1}=(R,S)_{2n+2}(R,S)_{2m+1}$\\ $=(R,S)_{2n+2m+3} $. When $R$ is preorder we  have: $(R,S)_{2n+1}(R,S)_{2m+1}$\\$=(R,S)_{2n}R^2(S,R)_{2m}=(R,S)_{2n}R(S,R)_{2m}=(R,S)_{2(n+m)+1}$. \\
The point 4) is straightforward consequence of the three first ones and the point 5) is obtained by straightforward inductions.
\endproof
So, we get the following corollary which emphasizes the importance of $(R,S)_3$ $=RSR$ already observed in Proposition \ref{n=1}:
\begin{coro}\label{cor3}
	Let $\EE$ be a regular category, and $(R,S)$ any pair of reflexive relations on $X$.  Then we get:
	$(R,S)_{2k+1}\subset (RSR)^k \subset (R,S)_{4k-1}$;\\
	and when $R$ is a preorder: $(RSR)^k=(R,S)_{2k+1}$.
\end{coro}

\subsection{Suprema of the intertwined chains}

In this section we shall be interested in the existence of suprema of chains of reflexive relations. Starting with a chain $(R_i)_{i \in \mathbb N}$ of reflexive relations on $X$, we shall denote by $\Sigma_{i \in \mathbb N}R_i$ its supremum in the fiber $Ref_X\EE$ of reflexive relations on $X$. Straightforward are the following observations:
\begin{lemma}
	Let $\EE$ be a category with suprema of chains of reflexive relations: 
	1) these suprema preserve the inclusion of chains;\\
    2) these suprema preserve the passage to the dual;\\
    3) for any reflexive relation $T$, we get: $\Sigma_{i \in \mathbb N}(R_i\cap T) \subset (\Sigma_{i \in \mathbb N}R_i)\cap T$;\\ 
    4) for any reflexive relation $T$, we get: $\Sigma_{i \in \mathbb N}(R_iT) \subset (\Sigma_{i \in \mathbb N}R_i)T$.\\
	The point 1) implies that when a morphism $f:X\to Y$ in $\EE$ determines a morphism of chains of reflexive relations $(R_i)_{i \in \mathbb N}\to (T_i)_{i \in \mathbb N}$, it induces a morphism $\Sigma_{i \in \mathbb N}R_i \to  \Sigma_{i \in \mathbb N}T_i$ in $Ref\EE$ above $f$.
\end{lemma}
 When any chain of reflexive relations admits a supremum in $\EE$, the intertwining relations implies: $\Sigma((R,S))=\Sigma ((S,R))$ and $\Sigma^{op}((R,S))=\Sigma ((R^{op},S^{op}))$. So, when $R$ and $S$ are symmetric, so is $\Sigma((R,S))$. 
\begin{theo}\label{supermain}
	Let $\EE$ be a regular category. When any chain of reflexive relations admits a supremum which, in addition, is stable under composition of reflexive relations, then the supremum $\Sigma((R,S))$ of the chain $((R,S))$ is a preorder. Accordingly,  when, moreover, $R$ and $S$ are symmetric, $\Sigma ((R,S))$ is an equivalence relation.
\end{theo}
\proof
The stability under composition by a reflexive relation $T$ means:\\ $\Sigma_{i \in \mathbb N}(R_iT)=(\Sigma_{i \in \mathbb N}R_i)T$. Accordingly, we get:\\
$\Sigma_{i \in \mathbb N}((R,S)_i\Sigma((R,S)))=\Sigma^2((R,S))$. On the other hand, for any $k\in \mathbb{N}$, we have:
$(R,S)_{k+i-1}\subset (R,S)_{k}(R,S)_{i} \subset (R,S)_{k+i+1}$; whence:\\ $\Sigma ((R,S))\subset (R,S)_k\Sigma ((R,S))\subset \Sigma ((R,S))$, and $(R,S)_k\Sigma ((R,S))=\Sigma ((R,S))$. Finally, we get: $\Sigma_{i \in \mathbb N}((R,S)_i\Sigma((R,S)))=\Sigma((R,S))$.
\endproof
From that we get a very useful construction of the preorder and the equivalence relation generated by a reflexive relation:
\begin{coro}
	Let $\EE$ be a regular category and suppose that any chain of reflexive relations admits a supremum which is stable under composition of reflexive relations. Then the inclusion $PrO\EE \hookrightarrow Ref\EE$ and $Equ\EE \hookrightarrow Ref\EE$ have left adjoints.
\end{coro}
\proof
The first adjoint associates the preorder $\Sigma ((T,T))$ with the reflexive relation $T$, while the second one associates the equivalence relation $\Sigma ((T,T^{op}))$ with it.
\endproof

\section{Suprema of preorders}

\subsection{Suprema of preorders and equivalence relations}

From Theorem \ref{supermain} we can immediately derive the following:

\begin{coro}\label{cofib}
	Let $\EE$ be a regular category and suppose that any chain of reflexive relations admits a supremum which is stable under composition of reflexive relations. Then, for any pair $((R,S))$ of preorders (resp. equivalence relations) on $X$, we get $\Sigma((R,S))=R\vee S$ in $PrO\EE$ (resp. in $Equ\EE$). So, any regular epimorphism $f:X\onto Y$ in $\EE$ admits a cocartesian map above it with any domain $T$ in $PrO\EE$ (resp. in $Equ\EE$). It is  given by: $T\hookrightarrow \Sigma ((R[f],T)) \twoheadrightarrow f(\Sigma ((R[f],T)))$.
\end{coro}
\proof
We know that $R\subset \Sigma ((R,S))$ and $S\subset \Sigma ((R,S))$. Given any preorder $T$ such that $R\subset T$ and $S\subset T$, we get $(R,S)_i\subset T$ for all $i\in \mathbb{N}$, so that: $\Sigma ((R,S)) \subset T$. Furthermore, when $R$ and $S$ are symmetric, we noticed that $\Sigma ((R,S))$ is an equivalence relation. The last assertion is then a consequence of Proposition \ref{prop10}.
\endproof

It remains to produce examples of regular category admitting suprema of chains  of reflexive relations which are stable under composition of reflexive relations.
\begin{enumerate}
	\item In the category $Set$ of sets, given any chain $(R_i)_{i\in \mathbb{N}}$ of reflexive relations, the union of the members of this chain is the supremum in question. It is easy to check that it is stable under composition of reflexive relations.
	\item Clearly, the same result holds in any variety $\VV$ of Universal Algebra. 
\end{enumerate}

\begin{coro}
	Given any variety $\VV$ of Universal Algebra, the fibration $(\;)_0: PrO\VV \to \VV$ and its restriction to $Equ\VV$ are  cofibrations as well.
\end{coro}
\proof
There are cocartesian maps above any regular epimorphism by Proposition \ref{cofib}; and above any monomorphism $u: U\into X$
as well, by the following construction: starting with a preorder (resp. equivalence  relation) $T$ on the algebra $U$, consider the family $\mathcal F$ of preorders (resp. equivalence  relations)  $S$ on $X$ such that $T\subset u^{-1}(S)$; then take the intersection $F$ of the preorders (resp. equivalence  relations) of this family. Then the induced map: $T \into F$ above $u$ is necessarily cocartesian.
\endproof
Here is the  main result of this section, which introduces a simple condition implying the assumption of Corollary \ref{cofib}:
\begin{prop}\label{chos}
	Let $\EE$ be a regular category admitting suprema of chains  of subobjects which are stable under pullback  along any split epimorphisms. Then $\EE$ has suprema of chains  of reflexive relations which are stable under composition  of reflexive relations. Acccordingly, $\EE$ has suprema of pairs of preorders and of equivalence relations, and the associated cocartesian maps.
\end{prop}
\proof
Any category $\EE$ admitting suprema of chains  of subobjects has  suprema of chains of reflexive relations. 

Now let $(S,T)$ be any pair of reflexive relations on $X$. Then the composition $ST$ of the reflexive relation $T$ with $S$ can be described in the following way:
$$ \xymatrix@=20pt{
	(d_0^S\times X)((d_1^S\times X)^{*}(T)) \ar@{ >->}[d] && (d_1^S\times X)^{*}(T)\ar@{ >->}[d] \ar@{->>}[rr] \ar@{->>}[ll]  && T \ar@{ >->}[d]^{(d_0^T,d_1^T)} \\
	X\times X && S\times X \ar@{->>}[rr]_{d_1^S\times X} \ar@{->>}[ll]^{d_0^S\times X} && X\times X
}
$$
with $ST=(d_0^S\times X)((d_1^S\times X)^{*}(T))$, where the right hand side square is a pullback along the split epimorphism $d_1^S\times X$ while  the left hand side one is  produced by the direct image along the split epimorphism $d_0^S\times X$. Now, the pullback functor $(d_1^S\times X)^{*}$ preserves the suprema of chain of subobjects by assumption, while, in any regular category, the direct image functor preserves the suprema of chain of subobjects. So, the suprema of chains of reflexive relations are stable under composition with reflexive relations.
\endproof

\begin{enumerate}
	\item Any topos  $\EE$ with suprema of chain of subobjects admits suprema of chain of reflexive relations which are stable under composition of reflexive relations.
	
	\item Any Grothendieck topos  admits suprema of chain of reflexive relations which are stable under composition of reflexive relations.
\end{enumerate}

The first point is straightforward consequence of Proposition \ref{chos}, since, given any topos $\EE$:\\
- the category $\EE$ is exact and consequently regular;\\
- any base change functor $f^*$, having a left adjoint $\pi_f$, preserves colimits and consequently preserves suprema of chain of subobjects.\\
The second point is the consequence of the following:
\begin{prop}
	Let $\EE$ be an elementary topos admitting suprema of chains  of subobjects. Given any topology $j$ in $\EE$, the topos $sh_j\EE$ of $j$-sheaves has suprema of chains  of subobjects. Acccordingly, this is the case for any Grothendieck topos.
\end{prop}
\proof
The topologies $j$ in the topos $\EE$ are in bijection with the left exact idempotent monads $(T,\eta, \mu)$ on $\EE$. The monad $T$ being left exact, it is clear that, given any chain $(A_i)_{i\in \mathbb{N}}$ of subobjects of $X$ in the sub-topos $\EE^T$ of $T$-algebras, the subobject $T(\cup_{i\in \mathbb{N}}A_i)$ of $X$ is the supremum of the chain $(A_{i\in \mathbb{N}})$  in $\EE^T$. The last assertion concerning the Grothendieck topoi is a straighforward consequence of Proposition \ref{presh} below concerning the presheaf categories.
\endproof

By adding some observation of \cite{topos}, we get:

\begin{coro}
	Let $\EE$ be an elementary topos with suprema of chains  of subobjects. Then the fibration $(\;)_0: PrO\EE \to \EE$ and its restriction to $Equ\EE$ are  cofibrations as well.
\end{coro}
\proof
There are cocartesian maps above any regular epimorphism by the previous observations and Corollary \ref{cofib}; and above any monomorphism as well by Corollary 3.3 in \cite{topos}:
let  $u:U\into X$ be any monomorphism in $\EE$ and $T$ any preorder (resp. equivalence relation) on $U$. Then the following upward pushout of $s_0^T$ along $u$:
$$ \xymatrix@=30pt{
	{T\;} \ar@{>->}[r]^{\tilde u} \ar@<-2ex>[d]_{d_0^T} \ar@<2ex>[d]^{d_1^T} & {S\;}  \ar@<-2ex>@{.>}[d]_{d_0^S}\ar@<2ex>@{.>}[d]^{d_1^S} \\
	{U\;} \ar@{>->}[r]_{u} \ar[u]|{s_0^T}  & X \ar[u]|{s_0^S} 
}
$$
produces, from $T$, a preorder (resp. equivalence relation) $S$ on $X$. So, by construction, this morphism of preorders (resp. equivalence relations) is cocartesian above $u$.
\endproof

\subsection{Stability properties of the condition of Proposition \ref{chos}} 

\begin{prop}\label{presh}
	Let $\EE$ be a category admitting suprema of chains  of subobjects which are stable under pullback along split epimorphisms. Then this property holds in any functor category $\mathcal{F}(\DD,\EE)$. When, in addition, $\EE$ is regular, so is $\mathcal{F}(\DD,\EE)$. Under these two conditions, the functor category $\mathcal{F}(\DD,\EE)$ has suprema of reflexive  relations which are stable under composition of reflexive relations. Accordingly any presheaf category $\mathcal{F}(\DD^{op},Set)$ satisfies this property
\end{prop}
\proof
In $\mathcal{F}(\DD,\EE)$, the limits and colimits are levelwise. So, given a chain of reflexive relations $(R_i)_{i\in \mathbb{N}}$ on the functor $F$ in $\mathcal{F}(\DD,\EE)$, we can construct a levelwise union defined by $(\cup_{i\in \mathbb{N}}R_i)(X)=\cup_{i\in \mathbb{N}}R_i(X)$ for all $X\in \DD$. In the proof of our first assertion, all the needed constructions being levewise, the proof is straighforward. When $\EE$ is a regular category, the regular epimorphisms in $\mathcal{F}(\DD,\EE)$, are the  levelwise regular epimorphic natural transformations. Accordingly $\mathcal{F}(\DD,\EE)$ is regular.
\endproof
\begin{prop}
	Let $\EE$ be a category admitting suprema of chain of subobjects which are stable under pullback along split epimorphisms. Then this property holds in any category $\EE^C$ of coalgebras  of a left exact comonad on $\EE$. When, in addition, $\EE$ is regular, so is $\EE^C$. Accordingly, under these two conditions the category $\EE^C$ has suprema of reflexive  relations which are stable under composition of reflexive relations.
\end{prop}
\proof
Let $(C,\epsilon,\gamma)$ be a left exact comonad on $\EE$ and $U: \EE^C \to \EE$ the left exact and right exact conservative forgetful functor from the category $\EE^C$ of $C$-coalgebras. So, the functor $U$ preserves and reflects the subobjects.

By the universal property of a supremum of a chain of subobjects on $X$ in $\EE$, the supremum of the image by $U$ of a chain of subobjects in $\EE^C$ is canonically endowed with a structure of coalgebras, so that $\EE^C$ has  suprema of chain of subobjects which are preserved and reflected by $U$. From that, our first assertion is straigtforward.

When $\EE$ is regular, a morphism $f$ in $\EE^C$ is regular if and only if its image $U(f)$ in $\EE$ is a regular one. From that the category $\EE^C$ of $C$-coalgebras is regular as well.  
\endproof

\section{Transfer of congruence modularity}

Let $\TT$ be the theory of any finitary algebraic structure, and $\VV(\TT)$ the associated variety. Let us denote by $\TT(\EE)$ the category of internal $\TT$-algebras in $\EE$ and  by $U_{\TT}^{\EE}:\TT(\EE)\to \EE$ the left exact and conservative forgetful functor associating with any $\TT$-algebra on an object $X\in \EE$, its ground object $X$. When  $\EE$ is regular and $f:X\onto Y$ is a regular epimorphism in $\EE$, so is $f^{\times n}: X^{\times n}\onto Y^{\times n}$. Accordingly, a morphism in $\TT(\EE)$ is a regular epimorphism if and only if its image by $U_{\TT}$ is a regular epimorphism. Therefore the category $\TT(\EE)$ is regular, and the functor $U_{\TT}$ preserves and reflects the direct images of subobjects; so it preserves and reflects the composition of relations.

In this section, we are going to give further precisions about the suprema $\Sigma((T,S))$ in $\TT(\EE)$ when $\EE$ is a regular category and  has suprema of chain of subobjects which are stable under pullback along  split epimorphisms and when, moreover, these suprema are stable \emph{under intersection with subobjects}. 
In particular, we shall show that, under this further condition, when the variety $\VV(\TT)$ is congruence modular, so is the category $\TT(\EE)$. This further condition clearly holds in any topos $\EE$ with suprema of chains of subobjects.

\begin{lemma}
	Let $\EE$ be a category with suprema of chains of subobjects which are stable under pullback along split epimorphisms and under intersection with subobjects. We get: $\Sigma_{i\in \NN}(A_i\times A_i)=\Sigma_{i\in \NN}(A_i)\times \Sigma_{i\in \NN}(A_i)$. So, when $p_n: X\times X \cdots X \to X$ is any $n$-ary operation on $X$, and $(A_i)_{i\in \mathbb{N}}$ a chain of subobjects on which this operation is stable, this operation is stable on the supremum $\Sigma_{i\in \mathbb{N}}((A_i))$ as well.
\end{lemma}
\proof
The stability of the suprema under intersection with subobjects means: $\Sigma_{i \in \mathbb N}(A_i\cap U)=(\Sigma_{i \in \mathbb N}A_i)\cap U$ for any family $(A_i,U)_{i\in \NN}$ of subobjects of $X$.

Starting with a chain $(A_i)_{i\in \mathbb{N}}$ of subobjects of $X$, we get: $\Sigma_{i\in \NN}(A_i\times A_i)\subset\Sigma_{i\in \NN}(A_i)\times \Sigma_{i\in \NN}(A_i)$ in any case. Now let u: $U\into X$ be any subobject of $X$ and let us consider the following succesive right hand side pullbacks:
$$ \xymatrix@=20pt{
	U\times A_i\ar@{ >->}[d] \ar@{=}[rr] && U\times A_i \ar@{ >->}[d]  \ar@{ >->}[rr]&& X\times A_i\ar@{ >->}[d] \ar@{->>}[rr]  && A_i \ar@{ >->}[d] \\
	X\times X && \ar@{ >->}[ll]^{u\times X} U\times X \ar@{ >->}[rr]_{u\times X} && X\times X \ar@{->>}[rr]_{p_1^X} && X
}
$$
Under our assumption, since $p_1^X$ is a split epimorphism and $u\times X$ a subobject, we get: $\Sigma_{i\in \NN}(U\times A_i)=U\times \Sigma_{i\in \NN}(A_i)$ as subobjects of $X\times X$.
Since $A_j\subset A_i$ and consequently $A_j\times  A_i\subset A_i\times A_i$ when $j\leq i$, we get:\\ $A_j\times \Sigma_{i\in \NN}(A_i)= \Sigma_{i\in \NN}(A_j\times A_i)\subset\Sigma_{i\in \NN}(A_i\times A_i)$.\\ Whence: $\Sigma_{j\in \NN}(A_j)\times \Sigma_{i\in \NN}(A_i)= \Sigma_{j\in \NN}(A_j\times \Sigma_{i\in \NN}(A_i))\subset\Sigma_{i\in \NN}(A_i\times A_i)$. The last assertion is then straighforward. 
\endproof

\begin{prop}\label{DDay}
	Let $\EE$ be a category with suprema of chains of subobjects which are stable under pullback along split epimorphisms and stable under  intersection with subobjects. Then 
	any category $\TT(\EE)$ has suprema of chain of subobjects which are stable under intersections with subobjects and which are preserved and thus reflected by the conservative forgetful functor $U_{\TT}^{\EE}:\TT(\EE)\to \EE$.
	
	When, moreover, $\EE$ is regular, we get: $(T\vee S)\cap R=
	\Sigma_{i \in \mathbb N}((T,S)_i\cap R)$ in $PrO\EE$ and in $Equ\EE$, and this property is transfered to any category $\TT(\EE)$. This is the case, in particular, for any topos $\EE$ with suprema of chain of subobjects.
\end{prop}
\proof
Under our assumptions, any chain $(A_i)_{i\in \NN}$ of subalgebras of an algebra $X$ in $\TT(\EE)$ is a subalgebra according to the previous lemma. So, $\TT(\EE)$ has suprema of chain of subobjects which are stable under intersections with subobjects and these suprema are preserved and thus reflected by the conservative forgetful functor $U_{\TT}^{\EE}$.

When, moreoever, $\EE$ is regular, we get: $(T\vee S)\cap R=
\Sigma_{i \in \mathbb N}((T,S)_i\cap R)$ in $PrO\EE$ and in $Equ\EE$ by Proposition \ref{chos} and the stability of suprema under intersection. Now, this property is transfered to $\TT(\EE)$ by the preservation and the reflection of the compositions of relations and of the suprema of chain of subobjects by the forgetful functor $U_{\TT}^{\EE}$.

Suppose now that $\EE$ is a topos with suprema of chain of subsobjects. We already noticed that it is a regular category. Let $(A_i)_{i \in \NN}$ be a chain of subobjects of $X$ and $U$ another subobject of $X$. The inclusion $A_i\cap U \into X$ can be described in the following  way where the right hand square is a pullback:
$$ \xymatrix@=20pt{
	U\cap A_i \ar@{ >->}[d] && U\cap A_i\ar@{ >->}[d] \ar@{>->}[rr] \ar@{=}[ll]  && A_i \ar@{ >->}[d]^{a_i} \\
	X && {\; U\;} \ar@{>->}[rr]_{u} \ar@{>->}[ll]^{u} && {\;X}
}
$$
The functor $u^*:\EE/X\to \EE/U$ admitting a right adjoint $\pi_u$ in a topos, it preserves the  suprema  of chains of subobjects. And the composition by the monomorphism $u$ obviously preserves these suprema in any category $\EE$.
\endproof
The assertion of the previous proposition concerning the equivalence relation $(T\vee S)\cap R$ in $\TT(\EE)$ was already observed in \cite{Day}  for any variety $\VV$, i.e. when $\EE$ is the topos $Set$.\\

\noindent\textbf{Yoneda embedding for internal $\TT$-algebras}

As everybody knows, any object $X$ in a category $\mathbb E$ produces a functor: $$Y(X)=Hom_{\mathbb E}(-,X):\mathbb E^{op}\to Set$$ This, in turn, produces a fully faithful functor:
$\;\;\mathbb E \to \mathcal F(\mathbb E^{op},Set)\;\;$
which is called the \emph{Yoneda embedding\index{embedding -Yoneda}}. It is left exact when, in addition, the category $\mathbb E$ is finitely complete, the left exactness property being a synthetic translation of the universal properties of the finite limits.

Then any theory $\TT$ of finitary algebraic structure produces a factorization $Y_{\mathbb T}$ making the following diagram commute and making it a pullback:
\[\xymatrix@C=2pc@R=2pc{  \mathbb T(\EE) \ar@{-->}[rr]^{Y_{\mathbb T}} \ar[d]_{U_{\mathbb T}^{\EE}} && \mathcal F(\mathbb E^{op},\VV(\mathbb T)) \ar[d]^{\mathcal F( \mathbb E^{op}, \mathcal U_{\mathbb T})} \\
	\mathbb E \ar[rr]_Y && \mathcal F(\mathbb E^{op},Set)}\]
where $U_{\mathbb T}^{\EE}$ and $\mathcal U_{\mathbb T}$ are the induced forgetful functors, which are both left exact and conservative, and which, accordingly, are faithful and reflect finite limits as well. Since the square is a pullback, \emph{the functor $Y_{\mathbb T}$ is fully faithful and left exact}. It is called the Yoneda embedding for internal $\TT$-algebras.

Suppose now that, in the category $\EE$, any pair of equivalence relations has a supremum.  It  is said to be \emph{congruence modular} if, given any triple $(R,S,T)$ of equivalence relations on $X$, then we get $(R\bigvee S)\cap T=R\bigvee (S\cap T)$ provided that we have $R\subset T$. We noticed that, when $\EE$ is regular, the existence of suprema of pairs of equivalence relations is equivalent to the existence of cocartesian maps. Let us recall from Proposition 2.7 in \cite{BGour} that the  congruence modular regular categories can be characterized in terms  of cocartesian maps as well:
\begin{prop}\label{mod}
	Suppose that the category $\EE$ is regular and that any pair $(R,S)$ of equivalence relations has a supremum $R\vee S$ among the equivalence relations, or, equivalently, that $Equ\EE$ has cocartesian maps above regular epimorphisms in $\mathbb E$. Then $\EE$ is congruence modular if and only if these cocartesian maps are stable under pullbacks along the maps in the fibers of the fibration $(\;)_0:Equ\EE\to \EE$ (i.e. the inclusions of equivalence relations).
\end{prop}

In \cite{Gu}, a congruence modular variety $\VV$ is characterized by the Shifting Principle in any algebra $A$ of this variety: given any triple $(T,S,R)$ where $T$ and $R$ are internal equivalence relations such that $T\subset R$ and $S$ a internal reflexive and symmetric relation such that $S\cap R\subset T$, then, for any quadruple $(x,x',t,t')$ of elements of the algebra $A$, the property described by the following diagram holds:
$$ \xymatrix@=15pt{
	x \ar[r]^{S} \ar@{.>}@(l,l)[d]_T \ar[d]^R & t \ar[d]^T \\
	x' \ar[r]_{S} & t' 
}
$$ 
where the existence of the dotted arrow is implied by the existence of the plain ones. Categorically speaking, this property is characterized by the fact that the following commutative square indexed by $1$ is a pullback inside the variety $\VV$:
$$ \xymatrix@C=3pc@R=2pc{
	T \square S \ar@{ >->}[d] \ar@<-1,ex>[r]_{\delta_0^T}\ar@<+1,ex>[r]^{\delta_1^T}
	& T \ar@{ >->}[d] \ar[l]\\
	R \square S \ar@<-1,ex>[r]_{\delta_0^R}\ar@<+1,ex>[r]^{\delta_1^R}
	& R  \ar[l]
}
$$
see \cite{BG}, \cite{BGum} and \cite{BGGGum}. Actually, since $S$ is symmetric, it is equivalent to the fact that the square indexed by $0$ is a pullback.
\begin{theo}
Let $\TT$ be the theory of any finitary algebraic structure such that $\VV(\TT)$ is congruence modular. When $\EE$ is regular and satisfies the condition of the previous proposition, the category $\TT(\EE)$ is congruence modular as well. This is the case, in particular, for any topos  $\EE$ with suprema of chain of subobjects.
\end{theo}
\proof
Since $\VV(\TT)$ satisfies the Shifting Principle, so does  $\mathcal F( \mathbb E^{op}, \VV(\TT))$.
Since the functor $Y_{\mathcal T}$ is fully faithful and left exact, it preserves and reflects any left exact limit. So, $\EE(\TT)$ satisfies the Shifting Principle as well. By the previous proposition, our assumptions imply that $(T\vee S)\cap R=
\Sigma_{i \in \mathbb N}((T,S)_i\cap R)$ in $\TT(\EE)$ when $\EE$ is regular. Now, we can mimick the proof of Lemma 3.2 for varieties in \cite{Gu}, namely we have to check that $(T,S)_i\cap R\subset T\vee (S\cap R)$ for any $2\leq i\in \NN$. We shall do it by induction. Since $\TT(\EE)$ is regular, it is enough to check it in $Set$ thanks to the embedding theorem for regular categories \cite{Barr}. 

Since $T\subset R$ and $S\cap R\subset T$, we get $(T,S)_2\cap R\subset (T,S\cap R)\subset T\vee (S\cap R)$. Suppose now that we have $(T,S)_i\cap R\subset T\vee (S\cap R)$. First notice that $(T,S)_i$ is symmetric since $T$ and  $S$ are equivalence relations, and that we have $(T,S)_{i+1}\subset(T,S)_i(S,T)$. Suppose $(x,z)\in (T,S)_{i+1}\cap R$; then there is some $t$  satisfying the situation described by the following left hand side diagram:
$$ \xymatrix@=15pt{
	x \ar[rr]^{(T,S)_{i}} \ar[d]_R  && t \ar@{=}[d]  &&&& x \ar[d]^R \ar@{.>}@(l,l)[d]_{T\vee (S\cap R)} \ar[rr]^{(T,S)_{i}}   && t \ar@{=}[d] \\
	z \ar[rr]_{(T,S)} && t &&&& z \ar[rr]_{(T,S)_{i}} && t
}
$$
So, we get the right hand side one. Now since $(T,S)_i\cap R\subset T\vee (S\cap R)\subset R$, the Shifting Principle implies the existence of the dotted arrow.
\endproof

\section{Stationary intertwined chains}

We shall be interested here in the situation where one of the intertwined chains is stationary from some level $n$, which will be another way to produce the supremum of the chains. This will allow us, among other things, to characterize those preorders $T$ on $X$ such that, for a given regular epimorphism $f:X\onto Y$, we get $f(T)^{n+1}=f(T)^n$.

\subsection{Asymmetry between the odd and even cases}

\begin{lemma}\label{keylemma}
	Let $\EE$ be a regular category, $(R,S)$ any pair of reflexive relations on $X$ and $2\leq n$. Suppose $(R,S)_{n+1}=(R,S)_n \; (*)$; we then have:\\ 
	1) $(S,R)_{n+2}=(S,R)_{n+1}$ and $(R,S)_{n+3}=(R,S)_{n+2}$;\\
	2) so, for all $1\leq i$: $(S,R)_{n+2i}=(S,R)_{n+2i-1}$, $(R,S)_{n+2i+1}=(R,S)_{n+2i}$;\\
	3) $(S,R)_n\subset (R,S)_n$ and $(R,S)_{n+1}\subset (S,R)_{n+1}$;\\
	4) so, for all $0\leq i$: $(S,R)_{n+2i}\subset (R,S)_{n+2i}$, $(R,S)_{n+2i+1}\subset(S,R)_{n+2i+1}$.
	Whence the following diagram:
	$$ \xymatrix@=20pt{
& \cdots  (R,S)_n \ar@{=}[r]& (R,S)_{n+1} \ar@{ >->}[r] \ar@{ >->}[d] & (R,S)_{n+2} \ar@{=}[r] & (R,S)_{n+3}\ar@{ >->}[d]  \cdots   \\
& \cdots (S,R)_n \ar@{ >->}[r] \ar@{ >->}[u] &  (S,R)_{n+1} \ar@{=}[r] &  (S,R)_{n+2} \ar@{ >->}[r] \ar@{ >->}[u] & (S,R)_{n+3}   \cdots  
	}
	$$
\end{lemma}
\proof
Suppose $(*)$. We then get 1) by: $(S,R)_{n+2}=S(R,S)_{n+1}=S(R,S)_n=(S,R)_{n+1}$; whence $(R,S)_{n+3}=R(S,R)_{n+2}=R(S,R)_{n+1}=(R,S)_{n+2}$.\\
Then 2) is obtained by a simple induction from 1).\\
On the other hand, we get 3) by i): $(S,R)_n\subset (R,S)_{n+1}= (R,S)_n$ from $(*)$\ and by ii): $(R,S)_{n+1}\subset (S,R)_{n+2}=(S,R)_{n+1}$ from 1).\\
Then 4) is obtained from 2) by a simple induction. 
\endproof

\begin{lemma}\label{keylemma2}
	Let $\EE$ be a regular category, $(R,S)$ any pair of reflexive relations on $X$ and $2\leq n$. The following conditions are equivalent:\\ 
	1) $(R,S)_{n+2}=(R,S)_{n}$;\\
	2) the chain $((R,S))$ is stationary from $n$;\\
	3) the reflexive relation $(R,S)_{n}$ is a preorder.
\end{lemma}
\proof 
Suppose 1); from  $(R,S)_{n+1}=(R,S)_{n}$, we get: $(R,S)_{n+2i+1}=(R,S)_{n+2i}$ for all $1\leq i$, and from $(R,S)_{n+2}=(R,S)_{n+1}$: $(R,S)_{n+2i+2}=(R,S)_{n+2i+1}$ for all $1\leq i$ by the previous lemma. So, we get $(R,S)_{n+k}=(R,S)_{n}$ for all $1\leq k$.\\
Suppose 2); then: $(R,S)_{n}(R,S)_{n}\subset (R,S)_{2n+1}=(R,S)_{n}$.\\
Suppose 3); then: $(R,S)_{n+2}=(R,S)(R,S)_n\subset (R,S)_n(R,S)_n\subset (R,S)_n$.  
\endproof
When $n$ is odd and $R$ a preorder we can reduce the previous conditions:
\begin{theo}\label{stat1}
	Let $\EE$ be a regular category and $(R,S)$ any pair of a preorder and a reflexive relation on $X$.  When $1\leq m$, the following conditions are equivalent:\\
	1) $(R,S)_{2m+2}=(R,S)_{2m+1}$;\\
	2) the chain $((R,S))$ is stationary from $2m+1$;\\
	3) the reflexive  relation $(R,S)_{2m+1}$ is a preorder.\\
	4) $(S,R)_{2m+2}=(R,S)_{2m+1}$;\\
	5) we get the mixed subpermutability condition: $(S,R)_{2m+1}\subset (R,S)_{2m+1}$.\\
	So, under any of these conditions, when $S$ is a preorder as well, we get:\\
	$(R,S)_{2m+1}=R\vee S$ among the preorders. This is true among the equivalence relations when so are both $R$ and $S$.
\end{theo}
\proof
Suppose 1). We get: $(R,S)_{2m+3}=(R,S)_{2m+2}R=(R,S)_{2m+1}R$. And since $R$ is a preorder, we get: $(R,S)_{2m+1}R=(R,S)_{2m+1}$. Whence: $(R,S)_{2m+1}$ $=(R,S)_{2m+2}=(R,S)_{2m+3}$. Now, by the previous lemma we get 2) and 3).

Suppose 3). Then: $(S,R)_{2m+2}\subset (R,S)_{2m+3}=(R,S)_{2m+1}(S,R)$;\\
and: $(R,S)_{2m+1}(S,R)\subset (R,S)_{2m+1}(R,S)_{2m+1}=(R,S)_{2m+1}\subset (S,R)_{2m+2}$.

Suppose 4). Then: $(S,R)_{2m+1}\subset (S,R)_{2m+2}= (R,S)_{2m+1}$.

Suppose 5). Then: $(R,S)_{2m+2}=R(S,R)_{2m+1}\subset R(R,S)_{2m+1}=(R,S)_{2m+1}$, since $R$ is preorder.

\smallskip

To prove the last sentence, let us show that, when $R$ and  $S$ are symmetric, so is $(R,S)_{2m+1}$. We get: 
$(R,S)_{2m+1}^{op}=(R^{op},S^{op})_{2m+1}=(R,S)_{2m+1}$.
\endproof

\begin{coro}\label{stat11}
	Let $\EE$ be a regular category, $f:X\onto Y$ a regular epimorphism and $T$ any reflexive relation on $X$. When $1\leq m$, the following conditions are equivalent:\\
	1) the reflexive relation $f(T)^m$ is a preorder;\\
	2) the pair $(R[f],T)$ is such that $(T,R[f])_{2m+1}\subset (R[f],T)_{2m+1}$.\\
	Under any of these conditions, when $T$ is a preorder (resp. an equivalence relation) the map:$\;\; T\hookrightarrow  (R[f],T)_{2m+1}\onto f((R[f],T)_{2m+1})$ $=f(T)^m$  is cocartesian above $f$ in $PrO\EE$ (resp. in $Equ\EE$).
\end{coro}
\proof
Since $R[f]$ is an equivalence relation, by the previous theorem we get 2) if and only if $(R[f],T)_{2m+1}$ is a preorder. Since this preorder contains $R[f]$, its direct image $f((R[f],T)_{2m+1})=f(T)^m$ (by Lemma \ref{image}) is a preorder by Theorem \ref{main3} .

Conversely, suppose that $f(T)^m$ is a preorder. By Proposition \ref{Rn}, it is equivalent to $f(T)^{m+1}=f(T)^m$. We get:\\ $(T,R[f])_{2m+1}\subset f^{-1}(f((T,R[f])_{2m+1}))= f^{-1}(f(T)^{m+1})=f^{-1}(f(T)^{m})$\\$=f^{-1}(f((R[f],T)_{2m+1}))= R[f](R[f],T)_{2m+1}R[f]=(R[f],T)_{2m+1}$\\
the penultimate equality being given  by Proposition \ref{idcomp}.\\
For the last assertion, apply Corollary \ref{firstex}.
\endproof

When $n$ is even, the situation is not so good:

\begin{theo}\label{stat2}
	Let $\EE$ be a regular category and $(R,S)$ any pair of a  reflexive relation and a preorder on $X$.  When $1\leq m$, the following conditions are equivalent:\\
	1) $(R,S)_{2m+1}=(R,S)_{2m}$;\\
	2) the chain $((R,S))$ is stationary from $2m$;\\
	3) the reflexive  relation $(R,S)_{2m}$ is a preorder.\\
	When $1\leq m$, the following conditions are equivalent:\\
	4) $(S,R)_{2m+1}=(R,S)_{2m}$;\\
	5) we get the mixed subpermutability condition: $(S,R)_{2m}\subset (R,S)_{2m}$.\\
	The second list is a consequence of the first  one.\\
	So, under any condition of the first list, when $R$ is a preorder as well, we get:
	$(R,S)_{2m}=R\vee S$ among the preorders. This is true among the equivalence relations, when so are both $R$ and $S$.
\end{theo}
\proof
Suppose 1). We get: $(R,S)_{2m+2}=(R,S)_{2m+1}S=(R,S)_{2m}S$. And since $S$ is a preorder, we get: $(R,S)_{2m}S=(R,S)_{2m}$. Whence: $(R,S)_{2m+2}$ $=(R,S)_{2m+1}=(R,S)_{2m}$. Now, by the Lemma \ref{keylemma2}, we get 2) and 3).

Suppose 3). Then: $(R,S)_{2m+1}\subset (R,S)_{2m+2}=(R,S)_{2m}(R,S)$;\\
$(R,S)_{2m}(R,S)\subset (R,S)_{2m}(R,S)_{2m}=(R,S)_{2m}\subset (R,S)_{2m+1}$.

Again suppose 1). Then: $(S,R)_{2m+1}\subset (R,S)_{2m+2}=(R,S)_{2m+1}S=(R,S)_{2m}.S=(R,S)_{2m}\subset (S,R)_{2m+1}$.

Suppose 4). Then: $(S,R)_{2m}\subset (S,R)_{2m+1}\subset (R,S)_{2m}$.

Suppose 5). Then: $(S,R)_{2m+1}=(S,R)_{2m}S\subset (R,S)_{2m}S=(R,S)_{2m}$, since $S$ is preorder.

\smallskip

To prove the last sentence, let us show that, when $R$ and  $S$ are symmetric, so is $(R,S)_{2m}$. We get: 
$(R,S)_{2m}^{op}=(S^{op},R^{op})_{2m}=(S,R)_{2m}\subset (R,S)_{2m}$.
\endproof

\begin{coro}
	Let $\EE$ be a regular category, $f:X\onto Y$ a regular epimorphism and $T$ any reflexive relation on $X$. When $1\leq m$ and $(T,R[f])_{2m+1}=(T,R[f])_{2m}\; (*)$,
	then the reflexive relation $f(T)^m$ is a preorder.
	
	Under this same condition, when $T$ is a preorder (resp. an equivalence relation) the map: $T\hookrightarrow  (T,R[f])_{2m}\onto f((T,R[f])_{2m})$ $=f(T)^m$  is cocartesian above $f$ in $PrO\EE$ (resp. in $Equ\EE$).
\end{coro}
\proof
Since $R[f]$ is an equivalence relation, by the previous theorem we get $(*)$ if and only if $(T,R[f])_{2m}$ is a preorder. Since this preorder contains $R[f]$, its direct image $f((T,R[f])_{2m})=f(T)^m$ (by Lemma \ref{image}) is a preorder by Theorem \ref{main3} .
For the last assertion, again apply Corollary \ref{firstex}.
\endproof

\subsection{Congruence $n$-permutable regular categories}

In this section, we shall globalize the results of the previous section, investigating another context where the existence of cocartesian maps is ensured: the $n$-permutable regular categories, namely the regular categories in which $(R,S)_n=(S,R)_n$ for any pair $(R,S)$ of equivalence relations. We shall incidentally  extend some results of \cite{ZRL} about this kind of categories.

\begin{theo}\label{disy}
	Let $\EE$ be a regular category and $1\leq m$. The following conditions are equivalent:\\
	1) $\EE$ is congruence $(2m+1)$-permutable;\\
	2) the direct image $f(R)$ along a regular epimorphism $f$ of any equivalence relation $R$ is such that $f(R)^m$ is an equivalence relation.
	
	When $\EE$ is congruence $2m$-permutable, the direct image $f(R)$ along a regular epimorphism $f$ of any equivalence relation $R$ is such that $f(R)^m$ is an equivalence relation.
	
	Accordingly, in any $n$-permutable category $\EE$, there are cocartesian images of equivalence relations.
	\end{theo}
\proof
For any $n\in \mathbb{N}$, the implication $1) \Rightarrow 2)$ is a straightforward consequence of Corollary \ref{stat11}. Conversely suppose 2). Given any reflexive relation $R$ on $X$, the direct image $d_1^R((d_0^R)^{-1}(S))$ of any reflexive relation $S$ on $X$ is nothing but $R\circ S\circ R^{op}$. So, when $R$ is an equivalence relation, it is $R\circ S\circ R$. In presence of 2),  $(d_1^R((d_0^R)^{-1}(S))^m=(R\circ S\circ R)^n=(R,S)_{2m+1}$ is an equivalence relation when so is $S$. Whence $(R,S)_{2m+1}=R\vee S$ for any pair of equivalence relations, and therefore $(R,S)_{2m+1}=(S,R)_{2m+1}$. Accordingly  the regular category $\mathbb E$ is congruence $(2n+1)$-permutable.
\endproof

Let us recall from \cite{ZRL} the following characterization theorem for the regular categories, see also the related result for varieties, announced in \cite{HM}:

\begin{theo}\label{zurab}
	Let $\EE$ be a regular category and  $3\leq  n$. The following conditions are equivalent:\\
	1) $\EE$ is congruence $n$-permutable;\\
	2) $(T,T^{op})_{n+1}=(T,T^{op})_{n-1}$ for any endorelation $T$;\\
	3) $R^{op}\subset R^{n-1}$ for any reflexive relation $R$;\\
	4) $R^n=R^{n-1}$ for any reflexive relation $R$.
\end{theo}
From 3), it is clear that, \emph{in a congruence $n$-permutable regular  category, any preorder is an equivalence relation}; so that $PrO\EE=Equ\EE$, and the questions of the existence of cocartesian images of preorders and equivalence relations coincide. From that we can produce another characterization which generalizes the well known one for Mal'tsev categories (case $n=1$), see \cite{CLP}:
\begin{coro}\label{genmal}
Let $\EE$ be a regular category. When $1\leq n$, the following conditions are equivalent:\\
1) $\EE$ is $n+1$-permutable;\\
2) any reflexive relation $R$ on $X$ is such that $R^{n}$ is an equivalence relation.
\end{coro}
\proof
The case $n=1$ is the Mal'tsev case. Suppose $2\leq n$.
By the point 4) of the previous theorem, by  Proposition \ref{Rn} and by our last observation, when $\EE$ is congruence $(n+1)$-permutable  then, for any reflexive relation $R$, the relation $R^{n}$ is an equivalence relation. Conversely, by  Proposition \ref{Rn}, the point 2) of our corollary implies that $R^{n+1}=R^n$. Accordingly, by the previous theorem, $\EE$ is congruence $(n+1)$-permutable.
\endproof
So, we get immediately:
\begin{prop}
	Let $\EE$ be a congruence $n$-permutable regular category. The inclusion functor: $Equ\EE \hookrightarrow Ref\EE$ has a left adjoint.
\end{prop}
\proof
Given any reflexive relation $T$ on $X$, then $T^{n-1}$ necessarily produces the associated equivalence relation.
\endproof

By the way, let us also recall here the following result from \cite{MRL} which can be applied to any congruence $n$-permutable regular category $\EE$:
\begin{prop}
	Let $\EE$ be a regular category. The following condition are equivalent:\\
	1) any preorder is an equivalence relation;\\
	2) any internal category is a groupoid.
\end{prop} 
In the case of varieties, according to \cite{Ch}, Condition 1) precisely characterizes those varieties $\VV$ which are $n$-permutable for some $2\leq n$. 

\subsection{Mixed $n$-subpermutability}

In Theorem 3.5 in \cite{CKP} is given a characterization of the congruence $n$-permutable regular categories  in terms of some  $(n-1)$-permutable reflexive relations. Trying to untangle the two different levels of permutation, the one of reflexive relations and the one of preorders, we get to another characterization which allows us to mix them:

\begin{theo}\label{mixed}
	Given any regular category $\EE$ and  any $2\leq n$, the following conditions are equivalent:\\
	1) the category $\EE$ is congruence $n$-permutable;\\
	2) the category $\EE$ satisfies the \emph{mixed $(2n-1)$-subpermutability condition}, i.e.:\\ $(S,R)_{2n-1}\subset (R,S)_{2n-1}$ for pairs $(R,S)$ of preorders and  reflexive relations.
\end{theo}
\proof
Suppose $n=2$; the mixed $3$-subpermutability is: $SRS\subset RSR$ for any pair for any pair $(R,S)$ of a preorder and a reflexive relation. Taking $R=\Delta_X$, we get $S^2\subset S$ for any reflexive relation, namely: any reflexive relation is a preorder. According to Proposition 1.2 in \cite{CPP}, this is equivalent to: any reflexive relation is an equivalence relation (i.e. $\EE$ is a Mal'tsev category), and when $\EE$ is regular, this is equivalent to the congruence $2$-permutability.

Conversely if $\EE$ is a Mal'tsev regular category, the reflexive relation $RSR=(R,S)_3$ is an equivalence relation which, by Theorem \ref{stat1}, is equivalent to the fact that $(S,R)_3\subset (R,S)_3$ when $R$ is a preorder.

Suppose 1) and $3\leq n$. Again taking $R=\Delta_X$, we get: $S^{n}\subset S^{n-1}$ for all reflexive relation $S$. 
By Theorem \ref{zurab}, we know that, in any regular category $\EE$, this last condition is equivalent to the  fact that the category $\EE$ is congruence $n$-permutable.\\
Now suppose 2) and $3\leq n$. By Corollary \ref{cor3} we know that, given any pair $(R,S)$ of a preorder and a reflexive relation, we have  $(R,S)_{2n+1}=(R,S)_3^{n}$. So, by the same Theorem \ref{zurab}, we get $(R,S)_{2n+1}=(R,S)_3^{n}=(R,S)_3^{n-1}=(R,S)_{2n-1}$. Now, by Theorem \ref{stat1}  since $R$ is a preorder, the category $\EE$ satisfies the mixed $(2n-1)$-subpermutability condition.
\endproof

\noindent keywords: internal reflexive relation, preorder and equivalence relation, supremum of pairs of internal preorders and of equivalence relations, cocartesian image, regular epimorphism and regular category, congruence modular variety and category, congruence $n$-permutable variety and category, elementary topos.\\
Mathematics Subject Classification: 18B25, 18C40, 18D30, 18E13, 08A30, 08B05.

\vspace{3mm}\noindent Univ. Littoral C\^ote d'Opale, UR 2597, LMPA,\\
Laboratoire de Math\'ematiques Pures et Appliqu\'ees Joseph Liouville,\\
F-62100 Calais, France.\\
bourn@univ-littoral.fr

\end{document}